\documentstyle[12pt]{article}

\input amssym.def
\input amssym
\newtheorem{definition}{Definition}
\newenvironment{remark}[1]{\vskip 5pt\noindent{\it Remark.\ }#1}{\vskip
5pt}

\newenvironment{proof}[1]{\noindent{\it Proof.\
}#1}{\hskip 3pt $\Box$\vskip 5pt}
\newtheorem{theorem}[definition]{Theorem}
\newtheorem{lemma}[definition]{Lemma}
\newtheorem{corollary}[definition]{Corollary}
\newtheorem{assertion}[definition]{Assertion}
\newtheorem{conjecture}[definition]{Conjecture}
\newtheorem{proposition}[definition]{Proposition}
\def\limfunc#1{\mbox{\rm #1}\,}
\def\text#1{\mbox{\rm #1}\,}
\def\stackunder#1#2{\mathrel{\mathop{#2}\limits_{#1}}}%

\begin{document}

\author{{\bf Steven Duplij} \thanks{%
Alexander von Humboldt Fellow} \thanks{%
On leave of absence from {\sl Theory Division, Nuclear Physics Laboratory,
Kharkov State University, KHARKOV 310077, Ukraine}} \thanks{%
E-mail: duplij@physik.uni-kl.de} \\
{\sl Physics Department, University of Kaiserslautern},\\
{\sl Postfach 3049, D-67653 KAISERSLAUTERN},\\
{\sl Germany}}
\title{{\bf SUPERMATRIX REPRESENTATIONS OF SEMIGROUP BANDS} }
\date{August 17, 1996 \hskip 1cm {\bf KL-TH-96/09}}
\maketitle

\begin{abstract}
Various semigroups of noninvertible supermatrices of the special
(antitriangle) shape having nilpotent Berezinian which appear in
supersymmetric theories are defined and investigated. A subset of them
continuously represents left and right zero semigroups and rectangular
bands. The ideal properties of higher order rectangular band
analogs and the ``wreath''
version of them are studied in detail. We introduce the ``fine'' equivalence
relations leading to ``multidimesional'' eggbox diagrams. They are full
images of Green's relations on corresponding subsemigroups.
\end{abstract}

\newpage

\section{Introduction}

Matrix semigroups \cite{pon,put1,putcha} are the great tool in concrete and
thorough investigation of detail abstract semigroup theory structure \cite
{cli/pre1,grillet,howie,ljapin}. Matrix representations
\cite{mca1,pra/sin} are
widely used in studying of finite semigroups \cite{lal/pet,rho/zal,zal},
topological semigroups \cite{bak/las,bro/fri2} and free semigroups \cite
{bren/char,fai1}. Usually matrix semigroups are defined over a field ${\Bbb K%
}$ \cite{okn1,pon1,pon2}. Nevertheless, after discovering of supersymmetry
\cite{vol/aku,wes/zum1} the realistic unified particle theories began to be
considered in superspace \cite{gat/gri/roc/sie,schw1}. In this picture all
variables and functions were defined not over a field ${\Bbb K}$, but over
Grassmann-Banach superalgebras over ${\Bbb K}$ \cite{rog1,dewitt} (or their
generalizations \cite{pes7,pes3,she}). However, the noninvertible (and
therefore semigroup) character of them was ignored for a long time, and only
recently the consistent studies of semigroups in supersymmetric theories
appeared \cite{dup6,dup11,dup10}. In addition to their physical contents
these investigations led to some nontrivial pure mathematical constructions
having unusual properties connected with noninvertibility and zero divisors
\cite{dup7,dup14}. In particular, it was shown \cite{dup12} that
supermatrices of the special shape can form
various strange and sandwich semigroups not
known before \cite{hic,mag/mis/tew}.

In this paper we work out continuous supermatrix representations of
semigroup bands\footnote{%
We note that study of idempotent semigroup representations \cite{siz},
especially the matrix ones \cite{erd}, is important by itself. The
idempotents also appear and are widely used in random matrix semigroup
applications \cite{ber2,dar/muk,muk1}.} introduced in \cite{dup12}. The
Green's relations on continuous zero semigroups and wreath rectangular bands
are studied in detail. We introduce the ``fine'' equivalence relations which
generalize them in some extent and lead to the ``multidimensional'' analog
of eggbox diagrams. Next investigations are connected with
continuous superanalogs of 0-simple semigroups and Rees's theorem which will
appear elsewhere.

\section{Preliminaries}

Let $\Lambda $ be a commutative Banach ${\Bbb Z}_2$-graded superalgebra over
a field ${\Bbb K}$ (where ${\Bbb K}={\Bbb R},$ ${\Bbb C}$ or ${\Bbb Q}_p$)
with a decomposition into the direct sum: $\Lambda =\Lambda _0\oplus \Lambda
_1$. The elements $a$ from $\Lambda _0$ and $\Lambda _1$ are homogeneous and
have the fixed even and odd parity defined as $\left| a\right| \stackrel{def%
}{=}\left\{ i\in \left\{ 0,1\right\} ={\Bbb Z}_2|\,a\in \Lambda _i\right\} $%
. If $\Lambda $ admit the decomposition into body and soul \cite{rog1} as $%
\Lambda ={\Bbb B}\oplus {\Bbb S}$, where ${\Bbb B}$ and ${\Bbb S}$ are
purely even and odd algebras over ${\Bbb K}$ respectively, the even
homomorphism ${\frak r}_{body}:\Lambda \rightarrow {\Bbb B}$ is called a
body map and the odd homomorphism ${\frak r}_{soul}:\Lambda \rightarrow
{\Bbb S}$ is called a soul map \cite{bar/bru/her,pes5}. Usually $\Lambda $
is modelled with the Grassmann algebras $\wedge \left( N\right) $ having $N$
generators \cite{rog1,vla/vol} or $\wedge \left( \infty \right) $ \cite
{boy/git,bry1,rog4}. The soul ${\Bbb S}$ is obviously a proper two-sided
ideal of $\Lambda $ which is generated by $\Lambda _1$. These facts allow us
to consider noninvertible morphisms on a par with invertible ones (in some
sense), which gives many interesting and nontrivial results (see e.g. \cite
{dup6,dup11,dup10}).

We consider $\left( p|q\right) $-dimensional linear model superspace $%
\Lambda ^{p|q}$ over $\Lambda $ (in the sense of \cite{berezin,leites}) as
the even sector of the direct product $\Lambda ^{p|q}=\Lambda _0^p\times
\Lambda _1^q$ \cite{rog1,vla/vol}. The even morphisms $\mbox{Hom}_0\left(
\Lambda ^{p|q},\Lambda ^{m|n}\right) $ between superlinear spaces $\Lambda
^{p|q}\rightarrow \Lambda ^{m|n}$ are described by means of $\left(
m+n\right) \times \left( p+q\right) $-supermatrices\footnote{The
supermatrix theory per se has own problems \cite
{bac/fel2,kob/nag1,hus/nie}, unexpected conclusions \cite{bac/fel1,rze1}
and renewed standard theorems \cite{urr/mor1}, which as
a whole attach
importance to more deep investigation of supermatrix systems
from the abstract viewpoint.}
(for details see \cite{berezin,lei1}).

\section{Supermatrix semigroups}

We consider $\left( 1+1\right) \times \left( 1+1\right) $-supermatrices
describing the elements from $\mbox{Hom}_0\left( \Lambda ^{1|1},\Lambda
^{1|1}\right) $ in the standard $\Lambda ^{1|1}$ basis \cite{berezin}
\begin{equation}
M\equiv \left(
\begin{array}{cc}
a & \alpha \\
\beta & b
\end{array}
\right) \in \mbox{{\rm Mat}}_\Lambda \left( 1|1\right)	\label{1}
\end{equation}

\noindent where $a,b\in \Lambda _0,\,\alpha ,\beta \in \Lambda _1$ (in the
following we use Latin letters for elements from $\Lambda _0$ and Greek
letters for ones from $\Lambda _1$, and all odd elements are nilpotent of
index 2). For sets of matrices we also use corresponding bold symbols, e.g. $%
{\bf M}\stackrel{def}{=}\left\{ M\in \mbox{{\rm Mat}}_\Lambda \left(
1|1\right) \right\} $, and the set product is ${\bf M\cdot N}\stackrel{def}{=%
}\left\{ \cup MN\,|\,M,N\in \mbox{{\rm Mat}}_\Lambda \left( 1|1\right)
\right\} $.

In this $\left( 1|1\right) $ case the supertrace defined as $\mbox{str}:%
\mbox{{\rm Mat}}_\Lambda \left( 1|1\right) \rightarrow \Lambda _0$ and
Berezinian (superdeterminant) defined as $\mbox{Ber}:\mbox{{\rm Mat}}%
_\Lambda \left( 1|1\right) \setminus \left\{ M|\,{\frak r}_{body}\left(
b\right) =0\right\} \rightarrow \Lambda _0$ are
\begin{equation}
\mbox{str}M=a-b,  \label{2}
\end{equation}
\begin{equation}
\mbox{Ber}M=\frac ab+\frac{\beta \alpha }{b^2}.  \label{3}
\end{equation}

In \cite{dup12} we introduced two kinds of possible reductions of $M$.

\begin{definition}
{\sl Even-reduced supermatrices} are elements from $\mbox{{\rm Mat}}_\Lambda
\left( 1|1\right) $ of the form
\begin{equation}
M_{even}\equiv \left(
\begin{array}{cc}
a & \alpha  \\
0 & b
\end{array}
\right) \in \mbox{{\rm RMat}}_\Lambda ^{\,even}\left( 1|1\right) \subset %
\mbox{{\rm Mat}}_\Lambda \left( 1|1\right) .  \label{4}
\end{equation}

\noindent {\sl Odd-reduced supermatrices} are elements from $\mbox{{\rm Mat}}%
_\Lambda \left( 1|1\right) $ of the form
\begin{equation}
M_{odd}\equiv \left(
\begin{array}{cc}
0 & \alpha  \\
\beta  & b
\end{array}
\right) \in \mbox{{\rm RMat}}_\Lambda ^{\,odd}\left( 1|1\right) \subset %
\mbox{{\rm Mat}}_\Lambda \left( 1|1\right) .  \label{5}
\end{equation}
\end{definition}

The odd-reduced supermatrices have a nilpotent Berezinian
\begin{equation}
\mbox{Ber}M_{odd}=\frac{\beta \alpha }{b^2}\Rightarrow \left( \mbox{Ber}%
M_{odd}\right) ^2=0  \label{5b}
\end{equation}

and satisfy
\begin{equation}
M_{odd}^n=b^{n-2}\left(
\begin{array}{cc}
\alpha \beta & \alpha b \\
\beta b & b^2-\left( n-1\right) \alpha \beta
\end{array}
\right) ,  \label{5c}
\end{equation}
which gives $\mbox{Ber}M_{odd}^n=0$ and $\limfunc{str}M_{odd}^n=b^{n-2}%
\left( n\alpha \beta -b^2\right) $.

It is seen that ${\bf M}$ is a set sum of ${\bf M}_{even}$ and ${\bf M}%
_{odd} $
\begin{equation}
{\bf M}={\bf M}_{even}{\bf \cup M}_{odd}{\bf .}  \label{6}
\end{equation}

The even- and odd-reduced supermatrices are mutually dual in the sense of
the Berezinian addition formula \cite{dup12}
\begin{equation}
\mbox{Ber}M=\mbox{Ber}M_{even}+\mbox{Ber}M_{odd}.  \label{7a}
\end{equation}

The matrices from $\mbox{Mat}\left( 1|1\right) $ form a linear semigroup of $%
\left( 1+1\right) \times \left( 1+1\right) $-supermatrices under the
standard supermatrix multiplication ${\frak M}\left( 1|1\right) \stackrel{def%
}{=}\left\{ {\bf M\,}|\,\cdot \right\} $ \cite{berezin}. Obviously, the
even-reduced matrices ${\bf M}_{even}$ form a semigroup ${\frak M}%
_{even}\left( 1|1\right) $ which is a subsemigroup of ${\frak M}\left(
1|1\right) $, because of ${\bf M}_{even}{\bf \cdot M}_{even}{\bf \subseteq M}%
_{even}$.

In general the odd-reduced matrices $M_{odd}$ do not form a semigroup, since
their multiplication is not closed
\begin{equation}
M_{odd\left( 1\right) }M_{odd\left( 2\right) }=\left(
\begin{array}{cc}
\alpha _1\beta _2 & \alpha _1b_2 \\
b_1\beta _2 & b_1b_2+\beta _1\alpha _2
\end{array}
\right) \notin {\bf M}_{odd}.  \label{15a}
\end{equation}

Nevertheless, some subset of ${\bf M}_{odd}$ can form a semigroup. Indeed,
due to the existence of zero divisors in $\Lambda $, from (\ref{15a}) it
follows that
\begin{equation}
\begin{array}{ccc}
{\bf M}_{odd}\cdot {\bf M}_{odd}\cap {\bf M}_{odd}\neq \emptyset &
\Rightarrow & \alpha \beta =0.
\end{array}
\label{15b}
\end{equation}

\begin{proposition}
1) The subsets ${\bf M}_{odd}|_{\alpha \beta =0}\subset {\bf M}_{odd}$ of
the odd-reduced matrices satisfying $\alpha \beta =0$ form a subsemigroup of
${\frak M}\left( 1|1\right) $ under the standard supermatrix multiplication.

2) In this semigroup the subset of matrices with $\beta =0$ is a left ideal,
and one with $\alpha =0$ is a right ideal, the matrices with $b=0$ form a
two-sided ideal.
\end{proposition}
\begin{proof}
Directly follows from (\ref{15a}).
\end{proof}

\begin{definition}
{\sl An odd-reduced semigroup }${\frak M}_{odd}\left( 1|1\right) $ is a
subsemigroup of ${\frak M}\left( 1|1\right) $ formed by the odd-reduced
matrices ${\bf M}_{odd}$ satisfying $\alpha \beta =0$.
\end{definition}

\section{One parameter subsemigroups of odd-reduced semigroup}

Let us investigate one-parameter subsemigroups of ${\frak M}_{odd}\left(
1|1\right) $. The simplest one is a semigroup of antidiagonal nilpotent
supermatrices of the shape
\begin{equation}
Y_\alpha \left( t\right) \stackrel{def}{=}\left(
\begin{array}{cc}
0 & \alpha t \\
\alpha & 0
\end{array}
\right) .  \label{16y}
\end{equation}

Together with a null supermatrix
\begin{equation}
Z\stackrel{def}{=}\left(
\begin{array}{cc}
0 & 0 \\
0 & 0
\end{array}
\right) .  \label{16z}
\end{equation}
they form a continuous null semigroup ${\frak Z}_\alpha \left( 1|1\right) $
having the null multiplication
\begin{equation}
Y_\alpha \left( t\right) Y_\alpha \left( u\right) =Z  \label{16yyz}
\end{equation}
(cf. \cite{cli/pre1}).

\begin{assertion}
For any fixed $t\in \Lambda ^{1|0}$ the set $\left\{ Y_\alpha \left(
t\right) ,Z\right\} $ is a 0-minimal ideal in ${\frak Z}_\alpha \left(
1|1\right) $.
\end{assertion}

In search of nontrivial one parameter subsemigroups ${\frak M}_{odd}\left(
1|1\right) $ we consider the odd-reduced supermatrices of the following
shape
\begin{equation}
P_\alpha \left( t\right) \stackrel{def}{=}\left(
\begin{array}{cc}
0 & \alpha t \\
\alpha & 1
\end{array}
\right)  \label{16}
\end{equation}
where $t\in \Lambda ^{1|0}$ is an even parameter of the Grassmann algebra $%
\Lambda $ which ''numbers'' elements $P_\alpha \left( t\right) $ and $\alpha
\in \Lambda ^{0|1}$ is a fixed odd element of $\Lambda $ which ''numbers''
the sets $\stackunder{t}{\cup }P_\alpha \left( t\right) $.

Here we will study one-parameter subsemigroups in ${\frak M}_{odd}\left(
1|1\right) $ as abstract semigroups \cite{cli/pre1,ljapin}, but not as
semigroups of operators \cite{davies,hil/phi}, which will be done elsewhere.

First, we establish multiplication properties of $P_\alpha \left( t\right) $
supermatrices. From (\ref{15a}) and (\ref{16}) it is seen that
\begin{equation}
\left(
\begin{array}{cc}
0 & \alpha t \\
\alpha & 1
\end{array}
\right) \left(
\begin{array}{cc}
0 & \alpha u \\
\alpha & 1
\end{array}
\right) =\left(
\begin{array}{cc}
0 & \alpha t \\
\alpha & 1
\end{array}
\right)  \label{m1}
\end{equation}

\begin{corollary}
The multiplication (\ref{m1}) is associative and so the $P_\alpha \left(
t\right) $ supermatrices form a semigroup ${\bf P}_\alpha $.
\end{corollary}

\begin{corollary}
All $P_\alpha \left( t\right) $ supermatrices are idempotent
\begin{equation}
\left(
\begin{array}{cc}
0 & \alpha t \\
\alpha	& 1
\end{array}
\right) ^2=\left(
\begin{array}{cc}
0 & \alpha t \\
\alpha	& 1
\end{array}
\right) .  \label{m12}
\end{equation}
\end{corollary}

\begin{proposition}
If $P_\alpha \left( t\right) =P_\alpha \left( u\right) $, then
\begin{equation}
t-u=\limfunc{Ann}\alpha .  \label{ma}
\end{equation}
\end{proposition}
\begin{proof}
From the definition (\ref{16}) it follows that two $P_\alpha \left( t\right)
$ supermatrices are equal iff $\alpha t=\alpha u$, which gives (\ref{ma}).
\end{proof}

Similarly we can introduce idempotent $Q_\alpha \left( t\right) $
supermatrices of the shape
\begin{equation}
Q_\alpha \left( t\right) \stackrel{def}{=}\left(
\begin{array}{cc}
0 & \alpha \\
\alpha t & 1
\end{array}
\right)  \label{16a}
\end{equation}
which satisfy
\begin{equation}
\left(
\begin{array}{cc}
0 & \alpha \\
\alpha t & 1
\end{array}
\right) \left(
\begin{array}{cc}
0 & \alpha \\
\alpha u & 1
\end{array}
\right) =\left(
\begin{array}{cc}
0 & \alpha \\
\alpha u & 1
\end{array}
\right)  \label{m1a}
\end{equation}
and form a semigroup ${\bf Q}_\alpha $.

\begin{assertion}
The semigroups ${\bf P}_\alpha $ and ${\bf Q}_\alpha $ contain no two sided
zeros and identities.
\end{assertion}

\begin{assertion}
The semigroups ${\bf P}_\alpha $ and ${\bf Q}_\alpha $ are continuous
unions of one
element groups with the action (\ref{m12}).
\end{assertion}

The relations (\ref{m1})--(\ref{m1a}) and
\begin{equation}
\left(
\begin{array}{cc}
0 & \alpha t \\
\alpha & 1
\end{array}
\right) \left(
\begin{array}{cc}
0 & \alpha \\
\alpha u & 1
\end{array}
\right) =\left(
\begin{array}{cc}
0 & \alpha t \\
\alpha u & 1
\end{array}
\right) \stackrel{def}{=}F_{tu},  \label{m1b}
\end{equation}
\begin{equation}
\left(
\begin{array}{cc}
0 & \alpha \\
\alpha u & 1
\end{array}
\right) \left(
\begin{array}{cc}
0 & \alpha t \\
\alpha & 1
\end{array}
\right) =\left(
\begin{array}{cc}
0 & \alpha \\
\alpha & 1
\end{array}
\right) \stackrel{def}{=}E  \label{m2b}
\end{equation}
are important from the abstract viewpoint and will be exploited below.

\begin{remark}
In general the supermatrix multiplication is noncommutative, noninvertible,
but associative, therefore any objects admitting supermatrix
representation (with closed multiplication) are automatically semigroups.
\end{remark}

\section{Continuous supermatrix representation of zero semigroups}

\label{sec4}

Let we consider an abstract set ${\frak P}_\alpha $ which consists of
elements ${\bf p}_t\in {\frak P}_\alpha $ ($t\in \Lambda ^{1|0}$ is a
continuous parameter) satisfying the multiplication law
\begin{equation}
{\bf p}_t*{\bf p}_u={\bf p}_t.	\label{p1}
\end{equation}

\begin{assertion}
The multiplication (\ref{p1}) is associative and therefore the set ${\frak P}%
_\alpha $ form a semigroup ${\cal P}_\alpha \stackrel{def}{=}\left\{ {\frak P%
}_\alpha ;*\right\} $.
\end{assertion}

\begin{assertion}
The semigroup ${\cal P}_\alpha $ is isomorphic to the left zero semigroup
\cite{cli/pre1} in which every element is both a left zero and a right
identity.
\end{assertion}

\begin{proposition}
The semigroup ${\cal P}_\alpha $ is epimorphic to the semigroup ${\bf P}%
_\alpha $.
\end{proposition}
\begin{proof}
Comparing (\ref{m1}) and (\ref{p1}) we observe that the mapping $\varphi :%
{\cal P}_\alpha \rightarrow {\bf P}_\alpha $ is a homomorphism. It is seen
that two elements ${\bf p}_t$ and ${\bf p}_u$ satisfying (\ref{ma}) have the
same image, i.e.
\begin{equation}
\varphi \left( {\bf p}_t\right) =\varphi \left( {\bf p}_u\right)
\Leftrightarrow t-u=\limfunc{Ann}\alpha ,\;{\bf p}_t,{\bf p}_u\in {\frak P}%
_\alpha .  \label{ker}
\end{equation}
\end{proof}

\begin{definition}
The relation
\begin{equation}
\Delta _\alpha =\left\{ \left( {\bf p}_t,{\bf p}_u\right) \,|\,t-u=\limfunc{%
Ann}\alpha ,\,{\bf p}_t,{\bf p}_u\in {\frak P}_\alpha \right\} .  \label{equ}
\end{equation}
is called $\alpha ${\sl -equality relation}.
\end{definition}

\begin{remark}
If the superparameter $t$ and $\alpha $ take value in different Grassmann
algebras which contain no mutually annihilating elements except zero, then $%
\limfunc{Ann}\alpha =0$ and $\Delta _\alpha =\Delta $.
\end{remark}

In the most of statements here the $\alpha $-equality relation $\Delta
_\alpha $ substitutes formally the standard equality relation $\Delta $,
nevertheless the fact that $\Delta \neq \Delta _\alpha $ leads to some new
structures and results. Among latter the following

\begin{corollary}
$\limfunc{Ker}\varphi =\stackunder{t\in \limfunc{Ann}\alpha }{\cup }{\bf p}_t
$.
\end{corollary}

\begin{remark}
Outside $\limfunc{Ker}\varphi $ the semigroup ${\cal P}_\alpha $ is
continuous and supersmooth, which can be shown by means of standard methods
of superanalysis \cite{berezin,dewitt}.
\end{remark}

\begin{assertion}
The semigroup ${\cal P}_\alpha $ is not reductive and not cancellative,
since ${\bf p}*{\bf p}_t={\bf p}*{\bf p}_u\Rightarrow {\bf p}_t\Delta
_\alpha {\bf p}_u$, but not ${\bf p}_t={\bf p}_u$ for all ${\bf p\in }{\frak %
P}_\alpha $. Therefore, the supermatrix representation given by $\varphi $
is not faithful.
\end{assertion}

\begin{corollary}
If $t+\limfunc{Ann}\alpha \cap u+\limfunc{Ann}\alpha \neq \emptyset $, then $%
{\bf p}_t\Delta _\alpha {\bf p}_u$.
\end{corollary}

Similarly the semigroup ${\cal Q}_\alpha $ with the multiplication
\begin{equation}
{\bf q}_t*{\bf q}_u={\bf q}_u  \label{p1q}
\end{equation}
is isomorphic to the right zero semigroup in which every element is both a
right zero and a left identity and epimorphic to the semigroup ${\bf Q}%
_\alpha $.

\begin{definition}
The semigroups ${\cal P}_\alpha $ and ${\cal Q}_\alpha $ can be named {\sl %
``somewhere commu\-ta\-tive''}\footnote{%
By analogy with nowhere commutative rectangular bands \cite{cli/pre1}.} or
{\sl ``almost anticom\-mutative''}, since for them ${\bf p}_t*{\bf p}_u={\bf p}%
_u*{\bf p}_t$ or ${\bf q}_t*{\bf q}_u={\bf q}_u*{\bf q}_t$ gives $\alpha
t=\alpha u$ and $t=u+\limfunc{Ann}\alpha $.
\end{definition}

\begin{proposition}
The semigroups ${\cal P}_\alpha $ and ${\cal Q}_\alpha $ are regular, but
not inverse.
\end{proposition}
\begin{proof}
For any two elements ${\bf p}_t$ and ${\bf p}_u$ using (\ref{p1}) we have $%
{\bf p}_t*{\bf p}_u*{\bf p}_t=\left( {\bf p}_t*{\bf p}_u\right) *{\bf p}_t=%
{\bf p}_t*{\bf p}_t={\bf p}_t$. Similarly for ${\bf q}_t$ and ${\bf q}_u$.
Then ${\bf p}_t$ has at least one inverse element ${\bf p}_u*{\bf p}_t*{\bf p%
}_u={\bf p}_u$. But ${\bf p}_u$ is arbitrary, therefore in semigroups ${\cal %
P}_\alpha $ and ${\cal Q}_\alpha $ any two elements are inverse. However,
they are not inverse semigroups in which every element has a unique inverse
\cite{cli/pre1}.
\end{proof}

The ideal structure of ${\cal P}_\alpha $ and ${\cal Q}_\alpha $ differs
somehow from the one of the left and right zero semigroups.

\begin{proposition}
Each element from ${\cal P}_\alpha $ forms by itself a principal right
ideal, each element from ${\cal Q}_\alpha $ forms a principal left ideal,
and therefore every principal right and left ideals in ${\cal P}_\alpha $
and ${\cal Q}_\alpha $ respectively have an idempotent generator.
\end{proposition}
\begin{proof}
From (\ref{p1}) and (\ref{p1q}) it follows that ${\bf p}_t={\bf p}_t*{\frak P%
}_\alpha $ and ${\bf q}_u={\frak Q}_\alpha *{\bf q}_u$.
\end{proof}

\begin{proposition}
The semigroups ${\cal P}_\alpha $ and ${\cal Q}_\alpha $ are left and right
simple respectively.
\end{proposition}
\begin{proof}
It is seen from (\ref{p1}) and (\ref{p1q}) that ${\frak P}_\alpha ={\frak P}%
_\alpha *{\bf p}_t$ and ${\frak Q}_\alpha ={\bf q}_u*{\frak Q}_\alpha $.
\end{proof}

The Green's relations on the standard left zero semigroup are the following:
${\cal L}$-equivalence coincides with the universal relation, and ${\cal R}$%
-equivalence coincides with the equality relation \cite{cli/pre1}. In our
case the first statement is the same, but instead of the latter we have

\begin{theorem}
In ${\cal P}_\alpha $ and ${\cal Q}_\alpha $ respectively ${\cal R}$%
-equivalence and ${\cal L}$-equivalence coincide with the $\alpha $-equality
relation (\ref{equ}). \label{th}
\end{theorem}
\begin{proof}
Consider the ${\cal R}$-equivalence in ${\cal P}_\alpha $. The elements $%
{\bf p}_t,{\bf p}_u\in {\cal P}_\alpha $ are ${\cal R}$-equivalent iff ${\bf %
p}_t*{\frak P}_\alpha ={\bf p}_u*{\frak P}_\alpha $. Using (\ref{equ}) we
obtain ${\cal R}=\Delta _\alpha $.
\end{proof}

\section{Wreath rectangular band}

Now we unify ${\cal P}_\alpha $ and ${\cal Q}_\alpha $ semigroups in some
nontrivial semigroup. First we consider the unified set of elements ${\frak P%
}_\alpha \cup {\frak Q}_\alpha $ and study their multiplication properties.
Using (\ref{m1b}) and (\ref{m2b}) we notice that ${\frak P}_\alpha \cap
{\frak Q}_\alpha ={\bf e}$, where $\varphi \left( {\bf e}\right) =E$ from (%
\ref{m2b}), and therefore ${\bf e\Delta }_\alpha {\bf p}_{t=1}$ and ${\bf %
e\Delta }_\alpha {\bf q}_{t=1}$. So we are forced to distinguish the region $%
t=1+\limfunc{Ann}\alpha $ from other points in the parameter superspace $%
\Lambda ^{1|0}$, and in what follows for any indices of ${\bf p}_t$ and $%
{\bf q}_t$ we imply $t\neq 1+\limfunc{Ann}\alpha $.

\begin{assertion}
${\bf e}$ is the left zero and right identity for ${\bf p}_t$, and ${\bf e}$
is the right zero and left identity for ${\bf q}_u$, i.e. ${\bf e*p}_t={\bf e%
},$ ${\bf p}_t*{\bf e=p}_t$, and ${\bf q}_u*{\bf e=e},$ ${\bf e*q}_u={\bf q}%
_u$.
\end{assertion}

Using (\ref{m2b}) it is easily to check that ${\bf q}_u{\bf *p}_t={\bf e}$,
but the reverse product needs to consider additional elements which are not
included in ${\frak P}_\alpha \cup {\frak Q}_\alpha $. From (\ref{m1b}) we
derive that
\begin{equation}
{\bf r}_{tu}{\bf =p}_t*{\bf q}_u,  \label{pq=r}
\end{equation}
where $\varphi \left( {\bf r}_{tu}\right) =F_{tu}$.

Let ${\frak R}_\alpha \stackrel{def}{=}\stackunder{t,u\neq 1+\limfunc{Ann}%
\alpha }{\cup }{\bf r}_{tu}$.

\begin{definition}
{\sl A wreath rectangular band }${\cal W}_\alpha $ is a set of idempotent
elements ${\frak W}_\alpha ={\frak P}_\alpha \cup {\frak Q}_\alpha \cup
{\frak R}_\alpha $ with a $*$-product (\ref{p1}) and the following Cayley
table

\medskip\

\begin{center}
\begin{tabular}{||c||c|c|c|c|c|c|c|c|c||}
\hline\hline
$1\setminus 2$ & ${\bf e}$ & ${\bf p}_t$ & ${\bf p}_u$ & ${\bf q}_t$ & ${\bf %
q}_u$ & ${\bf r}_{tu}$ & ${\bf r}_{ut}$ & ${\bf r}_{tw}$ & ${\bf r}_{vw}$ \\
\hline\hline
${\bf e}$ & ${\bf e}$ & ${\bf e}$ & ${\bf e}$ & ${\bf q}_t$ & ${\bf q}_u$ & $%
{\bf q}_u$ & ${\bf q}_t$ & ${\bf q}_w$ & ${\bf q}_w$ \\ \hline
${\bf p}_t$ & ${\bf p}_t$ & ${\bf p}_t$ & ${\bf p}_t$ & ${\bf r}_{tt}$ & $%
{\bf r}_{tu}$ & ${\bf r}_{tu}$ & ${\bf r}_{tt}$ & ${\bf r}_{tw}$ & ${\bf r}%
_{tw}$ \\ \hline
${\bf p}_u$ & ${\bf p}_u$ & ${\bf p}_u$ & ${\bf p}_u$ & ${\bf r}_{ut}$ & $%
{\bf r}_{uu}$ & ${\bf r}_{uu}$ & ${\bf r}_{ut}$ & ${\bf r}_{uw}$ & ${\bf r}%
_{uw}$ \\ \hline
${\bf q}_t$ & ${\bf e}$ & ${\bf e}$ & ${\bf e}$ & ${\bf q}_t$ & ${\bf q}_u$
& ${\bf q}_u$ & ${\bf q}_t$ & ${\bf q}_w$ & ${\bf q}_w$ \\ \hline
${\bf q}_u$ & ${\bf e}$ & ${\bf e}$ & ${\bf e}$ & ${\bf q}_t$ & ${\bf q}_u$
& ${\bf q}_u$ & ${\bf q}_t$ & ${\bf q}_w$ & ${\bf q}_w$ \\ \hline
${\bf r}_{tu}$ & ${\bf p}_t$ & ${\bf p}_t$ & ${\bf p}_t$ & ${\bf r}_{tt}$ & $%
{\bf r}_{tu}$ & ${\bf r}_{tu}$ & ${\bf r}_{tt}$ & ${\bf r}_{tw}$ & ${\bf r}%
_{tw}$ \\ \hline
${\bf r}_{ut}$ & ${\bf p}_u$ & ${\bf p}_u$ & ${\bf p}_u$ & ${\bf r}_{ut}$ & $%
{\bf r}_{uu}$ & ${\bf r}_{uu}$ & ${\bf r}_{ut}$ & ${\bf r}_{uw}$ & ${\bf r}%
_{uw}$ \\ \hline
${\bf r}_{tw}$ & ${\bf p}_t$ & ${\bf p}_t$ & ${\bf p}_t$ & ${\bf r}_{tt}$ & $%
{\bf r}_{tu}$ & ${\bf r}_{tu}$ & ${\bf r}_{tt}$ & ${\bf r}_{tw}$ & ${\bf r}%
_{tw}$ \\ \hline
${\bf r}_{vw}$ & ${\bf p}_v$ & ${\bf p}_v$ & ${\bf p}_v$ & ${\bf r}_{vt}$ & $%
{\bf r}_{vu}$ & ${\bf r}_{vu}$ & ${\bf r}_{vt}$ & ${\bf r}_{vw}$ & ${\bf r}%
_{vw}$ \\ \hline\hline
\end{tabular}
\end{center}

\medskip\

which is associative\footnote{%
For convenience and clearness we display some additional relations.} (as it
should be).
\end{definition}

From the Cayley table we can observe the following continuous subsemigroups
in the wreath rectangular band:

\begin{itemize}
\item  ${\bf e}$ - one element ``near identity'' subsemigroup;

\item  ${\tilde{{\cal P}}}_\alpha =\left\{ \stackunder{t\neq 1+\limfunc{Ann}%
\alpha }{\cup }{\bf p}_t;*\right\} $ -- ``reduced'' left zero semigroup;

\item  ${\cal P}_\alpha =\left\{ \stackunder{t\neq 1+\limfunc{Ann}\alpha }{%
\cup }{\bf p}_t\cup {\bf e};*\right\} $ -- full left zero semigroup;

\item  ${\tilde{{\cal Q}}}_\alpha =\left\{ \stackunder{t\neq 1+\limfunc{Ann}%
\alpha }{\cup }{\bf q}_t;*\right\} $ -- ``reduced'' right zero semigroup;

\item  ${\cal Q}_\alpha =\left\{ \stackunder{t\neq 1+\limfunc{Ann}\alpha }{%
\cup }{\bf q}_t\cup {\bf e};*\right\} $ -- full right zero semigroup;

\item  ${\tilde{{\cal F}}}_\alpha ^{\left( 1|1\right) }=\left\{ \stackunder{%
t,u\neq 1+\limfunc{Ann}\alpha }{\cup }{\bf r}_{tu};*\right\} $ --
``reduced'' rectangular band;

\item  ${\cal F}_\alpha ^{\left( 1|1\right) }=\left\{ \stackunder{t,u\neq 1+%
\limfunc{Ann}\alpha }{\cup }{\bf r}_{tu}\cup {\bf e};*\right\} $ -- full
rectangular band;

\item  ${\cal V}_\alpha ^L=\left\{ \stackunder{t,u\neq 1+\limfunc{Ann}\alpha
}{\cup }{\bf r}_{tu}\cup {\bf p}_t;*\right\} $ -- ``mixed'' left rectangular
band;

\item  ${\cal V}_\alpha ^R=\left\{ \stackunder{t,u\neq 1+\limfunc{Ann}\alpha
}{\cup }{\bf r}_{tu}\cup {\bf q}_t;*\right\} $ -- ``mixed'' right
rectangular band.
\end{itemize}

Thus we obtained the continuous supermatrix representation
for the left and right zero
semigroups and constructed from it the rectangular band supermatrix
representation. It is well known that any rectangular band isomorphic to a
Cartesian product of the left and right zero semigroups \cite{howie,petrich1}%
. Here we derived manifestly that (see (\ref{pq=r})) and presented the
concrete construction (\ref{m1b}). In addition, we unified all the above
semigroups in one object, viz. a wreath rectangular band.

\section{Rectangular band continuous representation}

The rectangular band multiplication is presented in the right lower corner
of the Cayley table. Usually \cite{cli/pre1,howie} it is defined by one
relation
\begin{equation}
{\bf r}_{tu}*{\bf r}_{vw}={\bf r}_{tw}.  \label{rr=r}
\end{equation}

In our case the indices are even Grassmann parameters from $\Lambda ^{1|0}$.
As for zero semigroups that also leads to some special peculiarities in the
ideal structure of such bands. Another difference is the absence of the
condition $u=v$ which arises in some applications from the finite nature of
indices considered as numbers of corresponding rows and columns in element
matrices (see e.g. \cite{lallement}).

Let us consider the Green's relations on ${\cal F}_\alpha ^{\left(
1|1\right) }$.

\begin{proposition}
Any two elements in ${\cal F}_\alpha ^{\left( 1|1\right) }$ are ${\cal J}$-
and ${\cal D}$-equivalent. \label{fd}
\end{proposition}
\begin{proof}
From (\ref{rr=r}) we derive
\begin{eqnarray}
{\bf r}_{tu}*{\bf r}_{vw}*{\bf r}_{tu}={\bf r}_{tw}*{\bf r}_{tu}={\bf r}%
_{tu},	\label{rrr=r} \\
{\bf r}_{vw}*{\bf r}_{tu}*{\bf r}_{vw}={\bf r}_{vw}*{\bf r}_{tw}={\bf r}_{vw}
\nonumber
\end{eqnarray}
for any $t,u,v,w\in \Lambda ^{1|0}$. First we notice that these equalities
coincide with the definition of ${\cal J}$-classes \cite{cli/pre1},
therefore any two elements are ${\cal J}$-equivalent, and so ${\cal J}$
coincides with the universal relation on ${\cal F}_\alpha ^{\left(
1|1\right) }$. Next using (\ref{rrr=r})  we observe that always ${\bf r}_{tu}%
{\cal R}{\bf r}_{tu}*{\bf r}_{vw}$ and ${\bf r}_{tu}*{\bf r}_{vw}{\cal L}%
{\bf r}_{vw}$. Since ${\cal D}={\cal L}\circ {\cal R}={\cal R}\circ {\cal L}$
(see e.g. \cite{howie}), then ${\bf r}_{tu}{\cal D}{\bf r}_{vw}$.
\end{proof}

\begin{assertion}
Every ${\cal R}$-class $R_{{\bf r}_{tu}}$consists of elements ${\bf r}_{tu}$
which are $\Delta _\alpha $-equivalent by the first index, i.e. ${\bf r}_{tu}%
{\cal R}{\bf r}_{vw}\Leftrightarrow t-v=\limfunc{Ann}\alpha $, and every $%
{\cal L}$-class $L_{{\bf r}_{tu}}$consists of elements ${\bf r}_{tu}$ which
are $\Delta _\alpha $-equivalent by the second index, i.e. ${\bf r}_{tu}%
{\cal L}{\bf r}_{vw}\Leftrightarrow u-w=\limfunc{Ann}\alpha $.
\end{assertion}
\begin{proof}
That follows from (\ref{rrr=r}), the manifest rectangular band decomposition
(\ref{pq=r}) and Theorem \ref{th}.
\end{proof}

So that the intersection of ${\cal L}$- and ${\cal R}$-classes is nonempty.
For the ordinary rectangular bands every ${\cal H}$-class consists of a
single element \cite{cli/pre1,howie}. In our case the situation is more
complicated.

\begin{definition}
The relation
\begin{equation}
\Delta _\alpha ^{\left( 1|1\right) }=\left\{ \left( {\bf r}_{tu},{\bf r}%
_{vw}\right) \,|\,t-v=\limfunc{Ann}\alpha ,\,u-w=\limfunc{Ann}\alpha ,\,{\bf %
r}_{tu},{\bf r}_{vw}\in {\frak R}_\alpha \right\} .  \label{dequ}
\end{equation}
is called a{\sl \ double }$\alpha ${\sl -equality relation}.
\end{definition}

\begin{theorem}
Every ${\cal H}$-class of ${\cal F}_\alpha ^{\left( 1|1\right) }$ consists
of double $\Delta _\alpha ^{\left( 1|1\right) }$-equivalent elements
satisfying ${\bf r}_{tu}\Delta _\alpha ^{\left( 2\right) }{\bf r}_{vw}$, and
so ${\cal H}=\Delta _\alpha ^{\left( 1|1\right) }$.
\end{theorem}
\begin{proof}
From (\ref{rrr=r}) and the definitions (\ref{m1b}) it follows that the
intersection of ${\cal L}$- and ${\cal R}$-classes happens when $\alpha
t=\alpha v$ and $\alpha u=\alpha w$. That gives $t=v+\limfunc{Ann}\alpha
,\,u=w+\limfunc{Ann}\alpha $ which coincides with the double $\alpha $%
-equality relation (\ref{dequ}).
\end{proof}

Let us consider the mapping $\psi :{\cal F}_\alpha ^{\left( 1|1\right)
}\rightarrow {\cal F}_\alpha ^{\left( 1|1\right) }/{\cal R}\times {\cal F}%
_\alpha ^{\left( 1|1\right) }/{\cal L}$ which maps an element ${\bf r}_{tu}$
to its ${\cal R}$- and ${\cal L}$-classes by
\begin{equation}
\psi \left( {\bf r}_{tu}\right) =\left\{ R_{{\bf r}_{tu}},L_{{\bf r}%
_{tu}}\right\} .  \label{psi}
\end{equation}
In the standard case $\psi $ is a bijection \cite{howie}. Now we have

\begin{assertion}
The mapping $\psi $ is a surjection.
\end{assertion}
\begin{proof}
That follows from Theorem \ref{th} and the decomposition (\ref{pq=r}).
\end{proof}

Let the Cartesian product ${\cal F}_\alpha ^{\left( 1|1\right) }/{\cal R}%
\times {\cal F}_\alpha ^{\left( 1|1\right) }/{\cal L}$ is furnished with the
rectangular band $\star $-multiplication of its ${\cal R}$- and ${\cal L}$%
-classes analogous to (\ref{rr=r}), i.e.
\begin{equation}
\left\{ R_{{\bf r}_{tu}},L_{{\bf r}_{tu}}\right\} \star \left\{ R_{{\bf r}%
_{vw}},L_{{\bf r}_{vw}}\right\} =\left\{ R_{{\bf r}_{tu}},L_{{\bf r}%
_{vw}}\right\} .  \label{l-r}
\end{equation}

For the standard rectangular bands the mapping $\psi $ is an isomorphism
\cite{howie,pas}. In our case we have

\begin{theorem}
The mapping $\psi $ is an epimorphism.
\end{theorem}
\begin{proof}
First we observe from (\ref{rrr=r}) that
\begin{eqnarray}
R_{{\bf r}_{tu}*{\bf r}_{vw}}=R_{{\bf r}_{tu}},  \label{rl} \\
L_{{\bf r}_{tu}*{\bf r}_{vw}}=L_{{\bf r}_{vw}},  \nonumber
\end{eqnarray}
and so under the $\star $-multiplication (\ref{l-r}) the mapping $\psi $ is
a homomorphism, since
\begin{eqnarray}
\psi \left( {\bf r}_{tu}*{\bf r}_{vw}\right) =\left\{ R_{{\bf r}_{tu}*{\bf r}%
_{vw}},L_{{\bf r}_{tu}*{\bf r}_{vw}}\right\} =\left\{ R_{{\bf r}_{tu}},L_{%
{\bf r}_{vw}}\right\}	\label{l-psi} \\
=\left\{ R_{{\bf r}_{tu}},L_{{\bf r}_{tu}}\right\} \star \left\{ R_{{\bf r}%
_{vw}},L_{{\bf r}_{vw}}\right\} =\psi \left( {\bf r}_{tu}\right) *\psi
\left( {\bf r}_{vw}\right) .  \nonumber
\end{eqnarray}

Then a surjective homomorphism is an epimorphism (e.g. \cite{hig1,hig}).
\end{proof}

\section{Higher $\left( n|n\right) $-band continuous representations}

Almost all above results can be generalized for the higher rectangular $%
\left( n|n\right) $ bands containing $2n$ continuous even Grassmann
parameters. The corresponding matrix construction is
\begin{equation}
F_{t_1t_2\ldots t_n,u_1u_2\ldots u_n}\stackrel{def}{=}\left(
\begin{array}{cc}
0 & \alpha t_1\;\alpha t_2\;\ldots \;\alpha t_n \\
\begin{array}{c}
\alpha u_1 \\
\alpha u_2 \\
\vdots \\
\alpha u_n
\end{array}
& I\left( n\times n\right)
\end{array}
\right) \in \mbox{{\rm RMat}}_\Lambda ^{\,odd}\left( 1|n\right) ,
\label{rmat-n}
\end{equation}
where $t_1,t_2\ldots t_n,u_1,u_2\ldots u_n\in \Lambda ^{1|0}$ are even
parameters, $\alpha \in \Lambda ^{1|0}$, $I\left( n\times n\right) $ is the
unit matrix, and the matrix multiplication is
\begin{equation}
F_{t_1t_2\ldots t_n,u_1u_2\ldots u_n}F_{t_1^{\prime }t_2^{\prime }\ldots
t_n^{\prime },u_1^{\prime }u_2^{\prime }\ldots u_n^{\prime
}}=F_{t_1t_2\ldots t_n,u_1^{\prime }u_2^{\prime }\ldots u_n^{\prime }}.
\label{ff=f}
\end{equation}

Thus the idempotent supermatrices $F_{t_1t_2\ldots t_n,u_1u_2\ldots u_n}$
form a semigroup ${\bf F}_\alpha ^{\left( n|n\right) }$.

\begin{definition}
{\sl A higher }$\left( n|n\right) ${\sl -band} ${\cal F}_\alpha ^{\left(
n|n\right) }\ni {\bf f}_{t_1t_2\ldots t_n,u_1u_2\ldots u_n}$ is represented
by the supermatrices $\mbox{{\rm RMat}}_\Lambda ^{\,odd}\left( 1|n\right) $
of the form (\ref{rmat-n}).
\end{definition}

The results of the Section \ref{sec4} with some slight differences hold
valid for ${\cal F}_\alpha ^{\left( n|n\right) }$ as well.

\begin{definition}
In ${\cal F}_\alpha ^{\left( n|n\right) }$ the relation
\begin{equation}
\begin{array}{c}
\Delta _\alpha ^{\left( n|n\right) }\stackrel{def}{=}\{\left( {\bf f}%
_{t_1t_2\ldots t_n,u_1u_2\ldots u_n},{\bf f}_{t_1^{\prime }t_2^{\prime
}\ldots t_n^{\prime },u_1^{\prime }u_2^{\prime }\ldots u_n^{\prime }}\right)
\,|\,t_k-t_k^{\prime }=\limfunc{Ann}\alpha , \\
\,u_k-u_k^{\prime }=\limfunc{Ann}\alpha ,\;1\leq k\leq n,\;{\bf f}%
_{t_1t_2\ldots t_n,u_1u_2\ldots u_n},{\bf f}_{t_1^{\prime }t_2^{\prime
}\ldots t_n^{\prime },u_1^{\prime }u_2^{\prime }\ldots u_n^{\prime }}\in
{\cal F}_\alpha ^{\left( 2n\right) }\}
\end{array}
\label{n-aple}
\end{equation}
is called {\sl a }$\left( n|n\right) ${\sl -ple }$\alpha ${\sl -equality
relation}.
\end{definition}

The semigroup ${\cal F}_\alpha ^{\left( n|n\right) }$ is also epimorphic to $%
{\bf F}_\alpha $, and two $\Delta _\alpha ^{\left( n|n\right) }$-equivalent
elements of ${\cal F}_\alpha ^{\left( n|n\right) }$ have the same image.

Let us consider $\mbox{{\rm RMat}}_\Lambda ^{\,odd}\left( k|m\right) $
idempotent supermatrices of the shape
\begin{equation}
F_{TU}\stackrel{def}{=}\left(
\begin{array}{cc}
0 & \alpha T \\
\alpha U & I
\end{array}
\right) ,  \label{tti}
\end{equation}
where $T\left( k\times m\right) $ and $U\left( m\times k\right) $ are the
band even parameter ordinary matrices and $I\left( m\times m\right) $ is the
unit matrix. This band contains maximum $2km$ parameters from $\Lambda
^{1|0} $.

The multiplication is
\begin{equation}
\left(
\begin{array}{cc}
0 & \alpha T \\
\alpha U & I
\end{array}
\right) \left(
\begin{array}{cc}
0 & \alpha T^{\prime } \\
\alpha U^{\prime } & I
\end{array}
\right) =\left(
\begin{array}{cc}
0 & \alpha T \\
\alpha U^{\prime } & I
\end{array}
\right) ,  \label{tui}
\end{equation}
which coincides in block form with the rectangular band multiplication (\ref
{rr=r})
\begin{equation}
F_{TU}F_{T^{\prime }U^{\prime }}=F_{TU^{\prime }}.  \label{ff=f1}
\end{equation}

\begin{theorem}
If $n=km$ the representations given by (\ref{rmat-n}) and (\ref{tti}) are
isomorphic.
\end{theorem}
\begin{proof}
Since in (\ref{ff=f}) and (\ref{ff=f1}) there exist no multiplication
between parameters, and so the representations given by matrices (\ref
{rmat-n}) and (\ref{tti}) differ by permutation if $n=km$.
\end{proof}

\begin{corollary}
The supermatrices $\mbox{{\rm RMat}}_\Lambda ^{\,odd}\left( 1|n\right) $ of
the shape (\ref{rmat-n}) exhaust all possible $\left( n|n\right) $-band
continuous representations.
\end{corollary}

\begin{remark}
The supermatrices (\ref{rmat-n}) represent $\left( k|m\right) $-bands as
well, where $1\leq k\leq n,\;1\leq m\leq n$. In this situation $t_{k+1}=1+%
\limfunc{Ann}\alpha ,\ldots t_n=1+\limfunc{Ann}\alpha ,\,u_{m+1}=1+\limfunc{%
Ann}\alpha ,\ldots u_n=1+\limfunc{Ann}\alpha $. So the above isomorphism
takes place for different bands having the same number of parameters
Therefore, we will consider below mostly the full $\left( n|n\right) $%
-bands, implying that they contain all particular and reduced cases.
\end{remark}

\begin{remark}
For $k=0$ and $m=0$ they describe $m$-right zero semigroups ${\cal Q}_\alpha
^{\left( m\right) }$ and $k$-left zero semigroups ${\cal P}_\alpha ^{\left(
k\right) }$ respectively having the following multiplication laws (cf. (\ref
{p1}) and (\ref{p1q}))
\begin{equation}
\begin{array}{c}
{\bf q}_{u_1u_2\ldots u_m}*{\bf q}_{u_1^{\prime }u_2^{\prime }\ldots
u_m^{\prime }}{\bf =q}_{u_1^{\prime }u_2^{\prime }\ldots u_m^{\prime }}, \\
{\bf p}_{t_1t_2\ldots t_k}*{\bf p}_{t_1^{\prime }t_2^{\prime }\ldots
t_k^{\prime }}{\bf =p}_{t_1t_2\ldots t_k}.
\end{array}
\label{qqpp}
\end{equation}
\end{remark}

\begin{proposition}
The $m$-right zero semigroups ${\cal Q}_\alpha ^{\left( m\right) }$ and $k$%
-left zero semigroups ${\cal P}_\alpha ^{\left( k\right) }$ are irreducible
in the sense that they cannot be presented as a direct product of
``1-dimen\-sional'' right zero ${\cal Q}_\alpha $ and left zero ${\cal P}%
_\alpha $ semigroups respectively.
\end{proposition}
\begin{proof}
It follows directly from comparing of the structure of supermatrices (\ref
{16}), (\ref{16a}) and (\ref{rmat-n}).
\end{proof}

\begin{proposition}
For the purpose of constructing $\left( k|m\right) $-bands one cannot use
``1-dimensional'' right zero ${\cal Q}_\alpha $ and left zero ${\cal P}%
_\alpha $ semigroups, because they reduce it to the ordinary
``2-dimensional'' rectangular band.
\end{proposition}
\begin{proof}
Indeed let ${\bf \tilde{f}}_{t_1t_2\ldots t_k,u_1u_2\ldots u_m}{\bf =p}%
_{t_1}*{\bf p}_{t_2}\ldots *{\bf p}_{t_k}*{\bf q}_{u_1}*{\bf q}_{u_2}\ldots *%
{\bf q}_{u_m}$. Then using the Cayley table above we derive ${\bf \tilde{f}}%
_{t_1t_2\ldots t_k,u_1u_2\ldots u_m}={\bf p}_{t_1}*{\bf q}_{u_m}$ which
trivially coincides with (\ref{pq=r}). Thus any combination of elements from
``1-dimensional'' right zero and left zero semigroups will not lead to new
construction other than in the Cayley table.
\end{proof}

Instead we have the following decomposition of a $\left( k|m\right) $-band
into $k$-left zero semigroup ${\cal P}_\alpha ^{\left( k\right) }$ and $m$%
-right zero semigroups ${\cal Q}_\alpha ^{\left( m\right) }$ $k$-left zero
semigroup ${\cal P}_\alpha ^{\left( k\right) }$
\begin{equation}
{\bf f}_{t_1t_2\ldots t_k,u_1^{\prime }u_2^{\prime }\ldots u_m^{\prime }}%
{\bf =p}_{t_1t_2\ldots t_k}*{\bf q}_{u_1^{\prime }u_2^{\prime }\ldots
u_m^{\prime }}.  \label{fpq}
\end{equation}

Despite this formula is similar to (\ref{pq=r}), we stress that the
increasing of number of superparameters is not an artificial trick, but a
natural way of searching for new constructions leading to generalization of
Green's relations and fine ideal structure of $\left( n|n\right) $-bands,
which has no analogs in the standard approach \cite{cli/pre1,howie,petrich1}.

\section{Fine ideal structure of $\left( n|n\right) $-bands}

Let us consider the Green's relations for $\left( n|n\right) $-bands. We
will try to establish the supermatrix meaning of properties of ${\cal R},%
{\cal L},{\cal D},{\cal H}$-classes. It will allow us to define and study
new equivalences most naturally, as well as to clear the previous
constructions. For clarity we use $\left( 2|2\right) $-band representation,
and the extending all the results to $\left( n|n\right) $-bands can be
easily done without further detail explanations.

The exact shape of the $\left( 2|2\right) $-band ${\cal F}_\alpha ^{\left(
2|2\right) }\ni {\bf f}_{t_1t_2,u_1u_2}$ supermatrix representation is
\begin{equation}
F_{t_1t_2,u_1u_2}=\left(
\begin{array}{ccc}
0 & \alpha t_1 & \alpha t_2 \\
\alpha u_1 & 1 & 0 \\
\alpha u_2 & 0 & 1
\end{array}
\right) .  \label{f22}
\end{equation}

According to the definition of ${\cal R}$-classes \cite{cli/pre1},
two elements $F_{t_1t_2,u_1u_2}$ and $F_{t_1^{\prime }t_2^{\prime
},u_1^{\prime }u_2^{\prime }}$ are ${\cal R}$-equivalent iff there exist two
another elements $X_{x_1x_2,y_1y_2}$, $W_{v_1v_2,w_1w_2}$ such that $%
F_{t_1t_2,u_1u_2}X_{x_1x_2,y_1y_2}=F_{t_1^{\prime }t_2^{\prime },u_1^{\prime
}u_2^{\prime }}$ and $F_{t_1^{\prime }t_2^{\prime },u_1^{\prime }u_2^{\prime
}}W_{v_1v_2,w_1w_2}=F_{t_1t_2,u_1u_2}$ simultaneously. In manifest form
\begin{equation}
F_{t_1t_2,u_1u_2}X_{x_1x_2,y_1y_2}=\left(
\begin{array}{ccc}
0 & \alpha t_1 & \alpha t_2 \\
\alpha y_1 & 1 & 0 \\
\alpha y_2 & 0 & 1
\end{array}
\right) =\left(
\begin{array}{ccc}
0 & \alpha t_1^{\prime } & \alpha t_2^{\prime } \\
\alpha u_1^{\prime } & 1 & 0 \\
\alpha u_2^{\prime } & 0 & 1
\end{array}
\right)  \label{fxf}
\end{equation}

and
\begin{equation}
F_{t_1^{\prime }t_2^{\prime },u_1^{\prime }u_2^{\prime
}}W_{v_1v_2,w_1w_2}=\left(
\begin{array}{ccc}
0 & \alpha t_1^{\prime } & \alpha t_2^{\prime } \\
\alpha w_1 & 1 & 0 \\
\alpha w_2 & 0 & 1
\end{array}
\right) =\left(
\begin{array}{ccc}
0 & \alpha t_1 & \alpha t_2 \\
\alpha u_1 & 1 & 0 \\
\alpha u_2 & 0 & 1
\end{array}
\right) .  \label{fwf}
\end{equation}

To satisfy the last equalities in (\ref{fxf}) and (\ref{fwf}) we should
choose
\begin{equation}
\begin{array}{c}
\alpha y_1=\alpha u_1^{\prime },\,\alpha y_1=\alpha u_1^{\prime }, \\
\alpha w_1=\alpha u_1,\,\alpha w_2=\alpha u_2,
\end{array}
\label{axaw}
\end{equation}
and
\begin{equation}
\alpha t_1=\alpha t_1^{\prime },\,\alpha t_2=\alpha t_2^{\prime }.
\label{atat}
\end{equation}

Due to the arbitrariness of $X_{x_1x_2,y_1y_2}$ and $W_{v_1v_2,w_1w_2}$ the
first equalities (\ref{axaw}) can be always solved by parameter
choice. The second equations (\ref{atat}) are the definition of ${\cal R}$%
-class of $\left( 2|2\right) $-band in the supermatrix interpretation. Thus
we have the following general

\begin{definition}
The ${\cal R}$-classes of $\left( n|n\right) $-band consist of elements
having all (!) $\alpha t_k$ fixed, where $1\leq k\leq n$.
\end{definition}

As the dual counterpart we formulate

\begin{definition}
The ${\cal L}$-classes of $\left( n|n\right) $-band consist of elements
having all (!) $\alpha u_k$ fixed, where $1\leq k\leq n$.
\end{definition}

In such a picture it is obvious that the join of these relations ${\cal D}=%
{\cal R}\vee {\cal L}$ covers all possible elements, and therefore any two
elements in $\left( n|n\right) $-band are ${\cal D}$-equivalent (cf.
Proposition \ref{fd}) . The intersection of them ${\cal H}={\cal R}\cap
{\cal L}$ obviously consists of the elements with all (!) $\alpha t_k$ and $%
\alpha u_k$ fixed. Indeed there is here the source of the $\left( n|n\right)
$-ple $\alpha $-equality relation definition (\ref{n-aple}).

\begin{proposition}
In $\left( 2|2\right) $-band ${\cal J}$-relation coincides with the
universal relation.
\end{proposition}
\begin{proof}
Multiplying (\ref{fxf}) by $F_{t_1t_2,u_1u_2}$ from the right and by $%
X_{x_1x_2,y_1y_2}$ from the left we obtain
\begin{equation}
\begin{array}{c}
F_{t_1t_2,u_1u_2}X_{x_1x_2,y_1y_2}F_{t_1t_2,u_1u_2}=F_{t_1t_2,u_1u_2}, \\
X_{x_1x_2,y_1y_2}F_{t_1t_2,u_1u_2}X_{x_1x_2,y_1y_2}=X_{x_1x_2,y_1y_2}
\end{array}
\label{fxfx}
\end{equation}
for any $t_1,t_2,u_1,u_2,x_1,x_2,y_1,y_2\in \Lambda ^{\left( 1|0\right) }$,
which coincides with the definition of ${\cal J}$-relation. The
arbitrariness of $F_{t_1t_2,u_1u_2}$ and $X_{x_1x_2,y_1y_2}$ proves the
statement.
\end{proof}

Summing up the standard approach for $\left( 2|2\right) $-bands we have
\begin{equation}
{\bf f}_{t_1t_2,u_1u_2}{\cal R}{\bf f}_{t_1^{\prime }t_2^{\prime
},u_1^{\prime }u_2^{\prime }}\Leftrightarrow \left\{ \alpha t_1=\alpha
t_1^{\prime }\wedge \alpha t_2=\alpha t_2^{\prime }\right\} ,  \label{r-st}
\end{equation}

\begin{equation}
{\bf f}_{t_1t_2,u_1u_2}{\cal L}{\bf f}_{t_1^{\prime }t_2^{\prime
},u_1^{\prime }u_2^{\prime }}\Leftrightarrow \left\{ \alpha u_1=\alpha
u_1^{\prime }\wedge \alpha u_2=\alpha u_2^{\prime }\right\} ,  \label{l-st}
\end{equation}

\begin{equation}
{\bf f}_{t_1t_2,u_1u_2}{\cal D}{\bf f}_{t_1^{\prime }t_2^{\prime
},u_1^{\prime }u_2^{\prime }}\Leftrightarrow \left\{
\begin{array}{c}
\left( \alpha t_1=\alpha t_1^{\prime }\wedge \alpha t_2=\alpha t_2^{\prime
}\right) \vee \\
\left( \alpha u_1=\alpha u_1^{\prime }\wedge \alpha u_2=\alpha u_2^{\prime
}\right)
\end{array}
\right\} ,  \label{d-st}
\end{equation}

\begin{equation}
{\bf f}_{t_1t_2,u_1u_2}{\cal H}{\bf f}_{t_1^{\prime }t_2^{\prime
},u_1^{\prime }u_2^{\prime }}\Leftrightarrow \left\{
\begin{array}{c}
\left( \alpha t_1=\alpha t_1^{\prime }\wedge \alpha t_2=\alpha t_2^{\prime
}\right) \wedge \\
\left( \alpha u_1=\alpha u_1^{\prime }\wedge \alpha u_2=\alpha u_2^{\prime
}\right)
\end{array}
\right\} .  \label{h-st}
\end{equation}

Now we are ready to introduce the fine ideal structure and understand what
was missed by the standard approach. From (\ref{r-st}) and (\ref{l-st}) it
is seen that the separate four possibilities for the equations to satisfy
are not covered by the ordinary ${\cal R}$- and ${\cal L}$-equivalent
relations. It is clear , why we wrote above exclamation marks: these
statements will be revised. So we are forced to define more general
relations, we call them ``fine equivalent relations''. They are appropriate
to describe all possible classes of elements in $\left( n|n\right) $-bands
missed by the standard approach. First we define them as applied to our
particular case for clarity.

\begin{definition}
The {\sl fine} ${\cal R}^{\left( k\right) }$- and ${\cal L}^{\left( k\right)
}${\sl -relations} on the $\left( 2|2\right) $-band are defined by
\begin{equation}
{\bf f}_{t_1t_2,u_1u_2}{\cal R}^{\left( 1\right) }{\bf f}_{t_1^{\prime
}t_2^{\prime },u_1^{\prime }u_2^{\prime }}\Leftrightarrow \left\{ \alpha
t_1=\alpha t_1^{\prime }\right\} ,  \label{r-my1}
\end{equation}
\begin{equation}
{\bf f}_{t_1t_2,u_1u_2}{\cal R}^{\left( 2\right) }{\bf f}_{t_1^{\prime
}t_2^{\prime },u_1^{\prime }u_2^{\prime }}\Leftrightarrow \left\{ \alpha
t_2=\alpha t_2^{\prime }\right\} ,  \label{r-my2}
\end{equation}
\begin{equation}
{\bf f}_{t_1t_2,u_1u_2}{\cal L}^{\left( 1\right) }{\bf f}_{t_1^{\prime
}t_2^{\prime },u_1^{\prime }u_2^{\prime }}\Leftrightarrow \left\{ \alpha
u_1=\alpha u_1^{\prime }\right\} ,  \label{l-my1}
\end{equation}
\begin{equation}
{\bf f}_{t_1t_2,u_1u_2}{\cal L}^{\left( 2\right) }{\bf f}_{t_1^{\prime
}t_2^{\prime },u_1^{\prime }u_2^{\prime }}\Leftrightarrow \left\{ \alpha
u_2=\alpha u_2^{\prime }\right\} .  \label{l-my2}
\end{equation}
\end{definition}

\begin{proposition}
The fine ${\cal R}^{\left( k\right) }$- and ${\cal L}^{\left( k\right) }$%
-relations are equivalence relations.
\end{proposition}
\begin{proof}
Follows from the manifest form of the multiplication and (\ref{fxf}) and (%
\ref{fwf}).
\end{proof}

Therefore, they divide ${\cal F}_\alpha ^{\left( 2|2\right) }$ to four fine
equivalence classes ${\cal F}_\alpha ^{\left( 2|2\right) }/{\cal R}^{\left(
k\right) }$ and ${\cal F}_\alpha ^{\left( 2|2\right) }/{\cal L}^{\left(
k\right) }$ as follows
\begin{equation}
R_{{\bf f}}^{\left( 1\right) }=\left\{ {\bf f}_{t_1t_2,u_1u_2}\in {\cal F}%
_\alpha ^{\left( 2|2\right) }\,|\,\alpha t_1=const\right\} ,  \label{rf1}
\end{equation}
\begin{equation}
R_{{\bf f}}^{\left( 2\right) }=\left\{ {\bf f}_{t_1t_2,u_1u_2}\in {\cal F}%
_\alpha ^{\left( 2|2\right) }\,|\,\alpha t_2=const\right\} ,  \label{rf2}
\end{equation}
\begin{equation}
L_{{\bf f}}^{\left( 1\right) }=\left\{ {\bf f}_{t_1t_2,u_1u_2}\in {\cal F}%
_\alpha ^{\left( 2|2\right) }\,|\,\alpha u_1=const\right\} ,  \label{lf1}
\end{equation}
\begin{equation}
L_{{\bf f}}^{\left( 2\right) }=\left\{ {\bf f}_{t_1t_2,u_1u_2}\in {\cal F}%
_\alpha ^{\left( 2|2\right) }\,|\,\alpha u_2=const\right\} .  \label{lf2}
\end{equation}

For clearness we can present schematically
\[
\begin{array}{ccccccccc}
\; & \; & \; & \; &  &	&  & R_{{\bf f}}^{\left( 1\right) } & R_{{\bf f}%
}^{\left( 2\right) } \\
&  &  &  &  &  &  & \updownarrow & \updownarrow
\end{array}
\]
\begin{equation}
\begin{array}{cc}
&  \\
L_{{\bf f}}^{\left( 1\right) } & \leftrightarrow \\
L_{{\bf f}}^{\left( 2\right) } & \leftrightarrow
\end{array}
\left(
\begin{array}{ccc}
0 & \alpha t_1 & \alpha t_2 \\
\alpha u_1 & 1 & 0 \\
\alpha u_2 & 0 & 1
\end{array}
\right) ,  \label{view}
\end{equation}
where arrows show which element of the supermatrix is fixed according to a
given fine equivalence relation.

From them we can obtain all known relations
\begin{equation}
{\cal R}^{\left( 1\right) }\cap {\cal R}^{\left( 2\right) }={\cal R},
\label{rrr}
\end{equation}
\begin{equation}
{\cal L}^{\left( 1\right) }\cap {\cal L}^{\left( 2\right) }={\cal L},
\label{lll}
\end{equation}
and
\begin{equation}
\left( {\cal R}^{\left( 1\right) }\cap {\cal R}^{\left( 2\right) }\right)
\cap \left( {\cal L}^{\left( 1\right) }\cap {\cal L}^{\left( 2\right)
}\right) ={\cal H},  \label{rrllh}
\end{equation}
\begin{equation}
\left( {\cal R}^{\left( 1\right) }\cap {\cal R}^{\left( 2\right) }\right)
\vee \left( {\cal L}^{\left( 1\right) }\cap {\cal L}^{\left( 2\right)
}\right) ={\cal D.}  \label{rrlld}
\end{equation}

However there are many other possible ``mixed'' equivalences which can be
classified using the definitions
\begin{equation}
{\cal H}^{\left( i|j\right) }={\cal R}^{\left( i\right) }\cap {\cal L}%
^{\left( j\right) },  \label{hrl}
\end{equation}
\begin{equation}
{\cal D}^{\left( i|j\right) }={\cal R}^{\left( i\right) }\vee {\cal L}%
^{\left( j\right) },  \label{drl}
\end{equation}
\begin{equation}
{\cal H}^{\left( ij|k\right) }=\left( {\cal R}^{\left( i\right) }\cap {\cal R%
}^{\left( j\right) }\right) \cap {\cal L}^{\left( k\right) },  \label{hrrl}
\end{equation}
\begin{equation}
{\cal H}^{\left( i|kl\right) }={\cal R}^{\left( i\right) }\cap \left( {\cal L%
}^{\left( k\right) }\cap {\cal L}^{\left( l\right) }\right) ,  \label{hrll}
\end{equation}
\begin{equation}
{\cal D}^{\left( ij|k\right) }=\left( {\cal R}^{\left( i\right) }\cap {\cal R%
}^{\left( j\right) }\right) \vee {\cal L}^{\left( k\right) },  \label{drrl}
\end{equation}
\begin{equation}
{\cal D}^{\left( i|kl\right) }={\cal R}^{\left( i\right) }\vee \left( {\cal L%
}^{\left( k\right) }\cap {\cal L}^{\left( l\right) }\right) .  \label{drll}
\end{equation}

The graphic interpretation of the mixed equivalence relations is given by
the following diagram

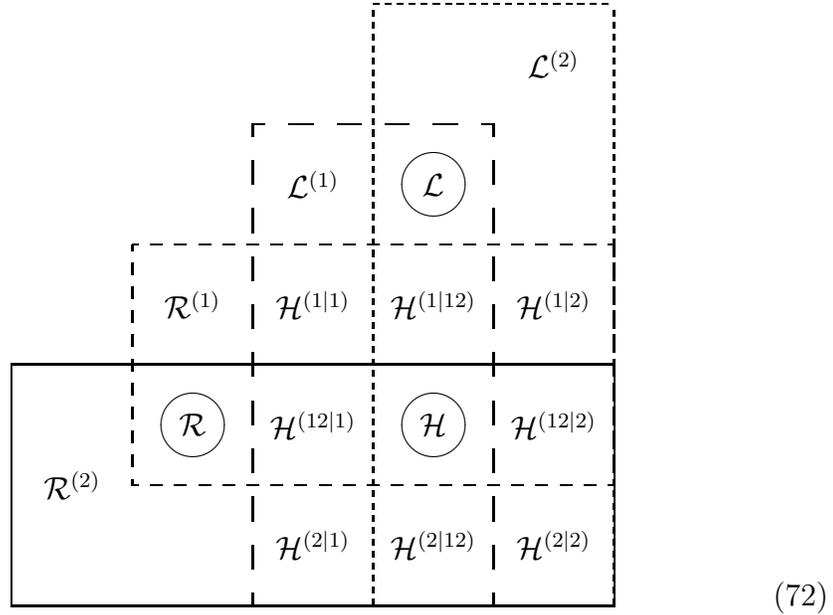
\begin{equation}
\setlength{\unitlength}{0.8cm}
\begin{picture}(10,10)
\put(0,0){\framebox(10,4){}}
\put(6,0){\dashbox{0.1}(4,10){}}
\put(2,2){\dashbox{0.2}(8,4){}}
\put(4,0){\dashbox{0.4}(4,8){}}
\put(6.5,2.5){\makebox(1,1){$\cal H$}}
\put(4.5,2.5){\makebox(1,1){${\cal H}^{\left( 12|1\right)}$}}
\put(4.5,0.5){\makebox(1,1){${\cal H}^{\left( 2|1\right)}$}}
\put(6.5,0.5){\makebox(1,1){${\cal H}^{\left( 2|12\right)}$}}
\put(8.5,0.5){\makebox(1,1){${\cal H}^{\left( 2|2\right)}$}}
\put(8.5,2.5){\makebox(1,1){${\cal H}^{\left( 12|2\right)}$}}
\put(8.5,4.5){\makebox(1,1){${\cal H}^{\left( 1|2\right)}$}}
\put(6.5,4.5){\makebox(1,1){${\cal H}^{\left( 1|12\right)}$}}
\put(4.5,4.5){\makebox(1,1){${\cal H}^{\left( 1|1\right)}$}}
\put(0.5,1.5){\makebox(1,1){${\cal R}^{\left( 2\right)}$}}
\put(2.5,4.5){\makebox(1,1){${\cal R}^{\left( 1\right)}$}}
\put(2.5,2.5){\makebox(1,1){${\cal R}$}}
\put(6.5,6.5){\makebox(1,1){${\cal L}$}}
\put(4.5,6.5){\makebox(1,1){${\cal L}^{\left( 1\right)}$}}
\put(8.5,8.5){\makebox(1,1){${\cal L}^{\left( 2\right)}$}}
\put(7,3){\circle{1}}
\put(7,7){\circle{1}}
\put(3,3){\circle{1}}
\end{picture}
\label{diag}
\end{equation}

\noindent
where the standard Green's relations are marked with circles. In (\ref{diag}%
) the standard ${\cal R}$- and ${\cal L}$-relations occupy 4
small squares
longwise, the ${\cal H}^{\left( i|j\right) }$-relations occupy 4
small squares in
square, the  ${\cal H}^{\left( ij|k\right) }$- and ${\cal H}^{\left(
i|jk\right) }$-relations occupy 2 small squares, the standard ${\cal H}$-relation
occupies 1 small square.

We observe that the mixed relations (\ref{hrl})-(\ref{drll}) are ``wider''
in some sense than the standard ones (\ref{rrr})-(\ref{rrlld}). Therefore,
using them we are able to describe thoroughly and appropriately all classes
of elements from $\left( n|n\right) $-bands including those which are missed
when one uses the standard Green's relations only\footnote{%
For nonnegative ordinary matrices the generalized Green's
relations (in some different sense) were
studied in \cite{yan/bar}.}.

For every mixed relation above we can determine a corresponding class using
obvious definitions. Then for every mixed ${\cal D}$-class we can build the
mixed eggbox diagram \cite{cli/pre1} of the fine ${\cal R,L}$-classes which
will have so many dimensions how many terms a given mixed relation has in
its right hand side of (\ref{drl}), (\ref{drrl}) and (\ref{drll}). For
instance, the eggbox diagrams of ${\cal D}^{\left( i|j\right) }$-classes are
two dimensional, but ones of ${\cal D}^{\left( ij|k\right) }$ and ${\cal D}%
^{\left( i|jk\right) }$-classes should be $3$-dimensional. In case of $%
\left( n|n\right) $-bands one has to consider all possible $k$-dimensional
eggbox diagrams, where $2\leq k\leq n-1$.

The introduced fine equivalence relations (\ref{r-my1})-(\ref{l-my2}) admit
a subsemigroup interpretation.

\begin{lemma}
The elements of ${\cal F}_\alpha ^{\left( n|n\right) }$ having $\alpha
t_k=\beta _k$ and $\alpha u_k=\gamma _k$, where $\beta _k,\gamma _k\in
\Lambda ^{0|1}$ are fixed, and $1\leq k\leq m$, form various $m$-index
subsemigroups.
\end{lemma}
\begin{proof}
It follows from the manifest form of matrix multiplication in (\ref{rmat-n}).
\end{proof}

We consider $\left( n-1\right) $-index subsemigroups of ${\cal F}_\alpha
^{\left( n|n\right) }$. They consist of elements having all but one $\alpha
t_k$ and all but one $\alpha u_k$ fixed. Let
\begin{equation}
{\cal U}_\alpha ^{\left( k\right) }\stackrel{def}{=}\left\{ {\bf f}%
_{t_1t_2\ldots t_n,u_1u_2\ldots u_n}\in {\cal F}_\alpha ^{\left( n|n\right)
}\,|\,\stackunder{i\neq k}{\wedge }\alpha t_i=\beta _i\stackunder{i\neq k}{%
\wedge }\alpha u_i=\gamma _i\right\}  \label{u}
\end{equation}
be a $\left( n-1\right) $-index subsemigroup which has only one nonfixed
pair $\alpha t_k$, $\alpha u_k$. The Green's relations on the subsemigroup $%
{\cal U}_\alpha ^{\left( k\right) }$ are the following
\begin{equation}
{\bf f}_{t_1t_2\ldots t_n,u_1u_2\ldots u_n}{\cal R}_U^{\left( k\right) }{\bf %
f}_{t_1^{\prime }t_2^{\prime }\ldots t_n^{\prime },u_1^{\prime }u_2^{\prime
}\ldots u_n^{\prime }}\Leftrightarrow \left\{ \alpha t_k=\alpha t_k^{\prime
}\right\} ,  \label{r-u}
\end{equation}
\begin{equation}
{\bf f}_{t_1t_2\ldots t_n,u_1u_2\ldots u_n}{\cal L}_U^{\left( k\right) }{\bf %
f}_{t_1^{\prime }t_2^{\prime }\ldots t_n^{\prime },u_1^{\prime }u_2^{\prime
}\ldots u_n^{\prime }}\Leftrightarrow \left\{ \alpha u_k=\alpha u_k^{\prime
}\right\} ,  \label{l-u}
\end{equation}
\begin{equation}
{\bf f}_{t_1t_2\ldots t_n,u_1u_2\ldots u_n}{\cal H}_U^{\left( k\right) }{\bf %
f}_{t_1^{\prime }t_2^{\prime }\ldots t_n^{\prime },u_1^{\prime }u_2^{\prime
}\ldots u_n^{\prime }}\Leftrightarrow \left\{ \alpha t_k=\alpha t_k^{\prime
}\wedge \alpha u_k=\alpha u_k^{\prime }\right\} ,  \label{h-u}
\end{equation}
\begin{equation}
{\bf f}_{t_1t_2\ldots t_n,u_1u_2\ldots u_n}{\cal D}_U^{\left( k\right) }{\bf %
f}_{t_1^{\prime }t_2^{\prime }\ldots t_n^{\prime },u_1^{\prime }u_2^{\prime
}\ldots u_n^{\prime }}\Leftrightarrow \left\{ \alpha t_k=\alpha t_k^{\prime
}\vee \alpha u_k=\alpha u_k^{\prime }\right\} ,  \label{d-u}
\end{equation}
where ${\bf f}_{t_1t_2\ldots t_n,u_1u_2\ldots u_n},{\bf f}_{t_1^{\prime
}t_2^{\prime }\ldots t_n^{\prime },u_1^{\prime }u_2^{\prime }\ldots
u_n^{\prime }}\in {\cal U}_\alpha ^{\left( k\right) }\subset {\cal F}_\alpha
^{\left( n|n\right) }$.

\begin{theorem}
The Green's relations on ${\cal U}_\alpha ^{\left( k\right) }$ are the
restrictions of the corresponding fine relations (\ref{r-my1})-(\ref{l-my2})
on ${\cal F}_\alpha ^{\left( n|n\right) }$ to the subsemigroup ${\cal U}%
_\alpha ^{\left( k\right) }$
\begin{equation}
{\cal R}_U^{\left( k\right) }={\cal R}^{\left( k\right) }\cap \left( {\cal U}%
_\alpha ^{\left( k\right) }\times {\cal U}_\alpha ^{\left( k\right) }\right)
,  \label{rruu}
\end{equation}
\begin{equation}
{\cal L}_U^{\left( k\right) }={\cal L}^{\left( k\right) }\cap \left( {\cal U}%
_\alpha ^{\left( k\right) }\times {\cal U}_\alpha ^{\left( k\right) }\right)
,  \label{lluu}
\end{equation}
\begin{equation}
{\cal H}_U^{\left( k\right) }={\cal H}^{\left( k|k\right) }\cap \left( {\cal %
U}_\alpha ^{\left( k\right) }\times {\cal U}_\alpha ^{\left( k\right)
}\right) ,  \label{hhuu}
\end{equation}
\begin{equation}
{\cal D}_U^{\left( k\right) }={\cal D}^{\left( k|k\right) }\cap \left( {\cal %
U}_\alpha ^{\left( k\right) }\times {\cal U}_\alpha ^{\left( k\right)
}\right) .  \label{dduu}
\end{equation}
\label{th-u}
\end{theorem}
\begin{proof}
It is sufficient to prove the statement for the particular case of ${\cal F}%
_\alpha ^{\left( 2|2\right) }$ and ${\cal U}_\alpha ^{\left( 1\right) }$,
and then to derive the general one by induction. Using the manifest form of $%
{\cal R}$-class definition (\ref{fxf})-(\ref{fwf}) we conclude that the
condition $\alpha t_1=\alpha t_1^{\prime }$ is common for the fine ${\cal R}%
^{\left( k\right) }$-classes and for the subsemigroup ${\cal R}_U^{\left(
k\right) }$-classes. By analogy one can prove other equalities.
\end{proof}

\begin{remark}
The second condition $\alpha t_2=\alpha t_2^{\prime }$ (which is a second
part of the definition of the ordinary ${\cal R}$-relation for ${\cal F}%
_\alpha ^{\left( 2|2\right) }$ (\ref{r-st})) holds in ${\cal U}_\alpha
^{\left( 1\right) }$ as well, but due to the subsemigroup own definition ($%
\alpha t_2=\beta _2=const$, $\alpha u_2=\gamma _2=const$), however $\alpha
t_2=\alpha t_2^{\prime }$ does not enter to the fine relation ${\cal R}%
^{\left( k\right) }$ at all. Therefore the latter is the most general one
among the ${\cal R}$-relations under consideration.
\end{remark}

\begin{remark}
The Theorem \ref{th-u} can be considered in view of \cite{hal3}, where the
formulas similar to (\ref{rruu})-(\ref{hhuu}) were proved, but with ordinary
Green's relations on the right hand side. Referring to the Diagram \ref{diag}
we conclude that our result contains the ordinary case \cite{hal3} as a
particular one.
\end{remark}

Moreover, we assume that the Theorem \ref{th-u} has more deep sense and
gives another treatment to the fine equivalence relations.

\begin{conjecture}
The Green's relations on a subsemigroup $U$ of $S$ have as counterpart
prolonged images in $S$ indeed the fine equivalence relations on $S$.\label
{conj}
\end{conjecture}

We proved this statement for the particular case of continuous $\left(
n|n\right) $-bands. It would be interesting to find and investigate other
possible algebraic systems where the Conjecture \ref{conj} is true.


\newpage

\begin{thebibliography}{10}

\bibitem{bac/fel2}
N.~B. Backhouse and A.~G. Fellouris,
\newblock {\em On the superdeterminant function for supermatrices},
\newblock J.~Phys.
\newblock {\bf  17} (1984), 1389--1395.

\bibitem{bac/fel1}
---------,
\newblock {\em Grassmann analogs of classical matrix groups},
\newblock J.~Math. Phys.
\newblock {\bf  26} (1985), 1146--1151.

\bibitem{bak/las}
J.~W. Baker and M.~Lashkarizadeh-{B}ami,
\newblock {\em On the representations of certain idempotent topological
  semigroups},
\newblock Semigroup Forum
\newblock {\bf  44} (1992), 245--254.

\bibitem{bar/bru/her}
C.~Bartocci, U.~Bruzzo, and D.~Hernandez Ruiperez,
\newblock {\em The Geometry of Supermanifolds},
\newblock Kluwer,
\newblock Dordrecht, 1991.

\bibitem{berezin}
F.~A. Berezin,
\newblock {\em Introduction to Superanalysis},
\newblock Reidel,
\newblock Dordrecht, 1987.

\bibitem{ber2}
M.~A. Berger,
\newblock {\em Central limit theorem for product of random matrices},
\newblock Trans. Amer. Math. Soc.
\newblock {\bf  285} (1984), 777--803.

\bibitem{boy/git}
C.~P. Boyer and S.~Gitler,
\newblock {\em The theory of ${G}^{\infty}$-supermanifolds},
\newblock Trans. Amer. Math. Soc.
\newblock {\bf  285} (1984), 241--267.

\bibitem{bren/char}
J.~L. Brenner and A.~Charnow,
\newblock {\em Free semigroups of $2 \times 2$ matrices},
\newblock Pacific J. Math.
\newblock {\bf  77} (1978), 57--69.

\bibitem{bro/fri2}
D.~R. Brown and M.~Friedberg,
\newblock {\em Linear representations of certain compact semigroups},
\newblock Trans. Amer. Math. Soc.
\newblock {\bf  160} (1971), 453--465.

\bibitem{bry1}
P.~Bryant,
\newblock {\em De {W}itt supermanifolds and infinite-dimensional ground rings},
\newblock J.~London Math. Soc.
\newblock {\bf  39} (1989), 347--368.

\bibitem{cli/pre1}
A.~H. Clifford and G.~B. Preston,
\newblock {\em The Algebraic Theory of Semigroups},
\newblock Vol.~1,
\newblock Amer. Math. Soc.,
\newblock Providence, 1961.

\bibitem{dar/muk}
R.~W.~R. Darling and A.~Mukherjea,
\newblock {\em Probability measures on semigroups of nonnegative matrices},
\newblock in The Analytical and Topological Theory of Semigroups,  (K.~H.
  Hofmann, J.~D. Lawson, and J.~S. Pym, eds.),
\newblock Walter de Gruyter, Berlin, 1990, pp.  361--377.

\bibitem{davies}
E.~B. Davies,
\newblock {\em One-Parameter Semigroups},
\newblock Academic Press,
\newblock London, 1980.

\bibitem{dup7}
S.~Duplij,
\newblock {\em On {$N=4$} super {R}iemann surfaces and superconformal
  semigroup},
\newblock J.~Phys.
\newblock {\bf  A24} (1991), 3167--3179.

\bibitem{dup6}
---------,
\newblock {\em On semigroup nature of superconformal symmetry},
\newblock J.~Math. Phys.
\newblock {\bf  32} (1991), 2959--2965.

\bibitem{dup11}
---------,
\newblock {\em Some abstract properties of semigroups appearing in
  superconformal theories},
\newblock University of Kaiserslautern, preprint KL-TH-95/11, hep-th-9505179
  (to appear in {\it Semigroup Forum}), 1995.

\bibitem{dup10}
---------,
\newblock {\em Ideal structure of superconformal semigroups},
\newblock Theor. Math. Phys.
\newblock {\bf  106} (1996), 355--374.

\bibitem{dup14}
---------,
\newblock {\em Noninvertibility and ''semi-'' analogs of (super) manifolds,
  fiber bundles and homotopies},
\newblock University of Kaiserslautern preprint KL-TH-96/10, 1996.

\bibitem{dup12}
---------,
\newblock {\em On an alternative supermatrix reduction},
\newblock Lett. Math. Phys.
\newblock {\bf  37} (1996), 385--396.

\bibitem{erd}
J.~A. Erdos,
\newblock {\em On products of idempotent matrices},
\newblock Glasgow Math. J.
\newblock {\bf  8} (1967), 118--122.

\bibitem{fai1}
V.~A. Faiziev,
\newblock {\em On pseudocharacters of a free semigroup invariant under its
  endomorphisms},
\newblock Russian Math. Surv.
\newblock {\bf  47} (1992), 205--206.

\bibitem{gat/gri/roc/sie}
S.~J. Gates, M.~T. Grisaru, M.~Rocek, et~al.,
\newblock {\em Superspace},
\newblock Benjamin,
\newblock Reading, 1983.

\bibitem{grillet}
P.-A. Grillet,
\newblock {\em Semigroups},
\newblock Marcel Dekker,
\newblock New York, 1994.

\bibitem{hal3}
T.E. Hall,
\newblock {\em Congruences and {G}reen's relations on regular semigroups},
\newblock Glasgow Math. J.
\newblock {\bf  11} (1972), 167--175.

\bibitem{hic}
J.~B. Hickey,
\newblock {\em Semigroups under sandwich operation},
\newblock Proc. Edinburgh Math. Soc.
\newblock {\bf  26} (1983), 371--382.

\bibitem{hig1}
P.~M. Higgins,
\newblock {\em A semigroup with an epimorphically embedded subband},
\newblock Bull. Amer. Math. Soc.
\newblock {\bf  27} (1983), 231--242.

\bibitem{hig}
---------,
\newblock {\em Completely semisimple semigroups and epimorphisms},
\newblock Proc. Amer. Math. Soc.
\newblock {\bf  96} (1986), 387--390.

\bibitem{hil/phi}
E.~Hille and R.~S. Phillips,
\newblock {\em Functional Analysis and Semigroups},
\newblock Amer. Math. Soc.,
\newblock Providence, 1957.

\bibitem{howie}
J.~M. Howie,
\newblock {\em An Introduction to Semigroup Theory},
\newblock Academic Press,
\newblock London, 1976.

\bibitem{hus/nie}
V.~Hussin and L.~M. Nieto,
\newblock {\em Supergroups factorizations through matrix realization},
\newblock J.~Math. Phys.
\newblock {\bf  34} (1993), 4199--4220.

\bibitem{kob/nag1}
Y.~Kobayashi and S.~Nagamishi,
\newblock {\em Characteristic functions and invariants of supermatrices},
\newblock J.~Math. Phys.
\newblock {\bf  31} (1990), 2726--2730.

\bibitem{lallement}
G.~Lallement,
\newblock {\em Semigroups and Combinatorial Applications},
\newblock Willey,
\newblock New York, 1979.

\bibitem{lal/pet}
G.~Lallement and M.~Petrich,
\newblock {\em Irreducible matrix representations of finite semigroups},
\newblock Trans. Amer. Math. Soc.
\newblock {\bf  139} (1969), 393--412.

\bibitem{leites}
D.~Leites,
\newblock {\em Supermanifold Theory},
\newblock Math. Methods Sci. Invest.,
\newblock Petrozavodsk, 1983.

\bibitem{lei1}
D.~A. Leites,
\newblock {\em Introduction to the theory of supermanifolds},
\newblock Russian Math. Surv.
\newblock {\bf  35} (1980), 1--64.

\bibitem{ljapin}
E.~S. Ljapin,
\newblock {\em Semigroups},
\newblock Amer. Math. Soc.,
\newblock Providence, 1968.

\bibitem{mag/mis/tew}
K.~D. Magill, P.~R. Misra, and U.~B. Tewari,
\newblock {\em Structure spaces for sandwich semigroups},
\newblock Pacific J. Math.
\newblock {\bf  99} (1982), 399--412.

\bibitem{mca1}
D.~B. McAlister,
\newblock {\em Representations of semigroups by linear transformations, 1,2},
\newblock Semigroup Forum
\newblock {\bf  2} (1971), 189--320.

\bibitem{muk1}
A.~Mukherjea,
\newblock {\em Convergence in distribution of products of random matrices: a
  semigroup approach},
\newblock Trans. Amer. Math. Soc.
\newblock {\bf  303} (1987), 395--411.

\bibitem{okn1}
J.~Okni\'{n}ski,
\newblock {\em Linear representations of semigroups},
\newblock in Monoids and Semigroups with Applications,  (J.~Rhodes, ed.),
\newblock World Sci., River Edge, 1991, pp.  257--277.

\bibitem{pas}
F.~Pastijn,
\newblock {\em Embedding semigroups in semibands},
\newblock Semigroup Forum
\newblock {\bf  14} (1977), 247--263.

\bibitem{pes7}
V.~Pestov,
\newblock {\em Ground algebras for superanalysis},
\newblock Rep. Math. Phys.
\newblock {\bf  29} (1991), 275--287.

\bibitem{pes3}
---------,
\newblock {\em Nonstandard hulls of normed {G}rassmannian algebras and their
  application in superanalysis},
\newblock Soviet Math. Dokl.
\newblock {\bf  317} (1991), 565--569.

\bibitem{pes5}
---------,
\newblock {\em Soul expansion of ${G}^{\infty}$ superfunctions},
\newblock J.~Math. Phys.
\newblock {\bf  34} (1993), 3316--3323.

\bibitem{petrich1}
M.~Petrich,
\newblock {\em Introduction to Semigroups},
\newblock Merill,
\newblock Columbus, 1973.

\bibitem{pon}
J.~S. Ponizovskii,
\newblock {\em On irreducible matrix semigroups},
\newblock Semigroup Forum
\newblock {\bf  24} (1982), 117--148.

\bibitem{pon2}
---------,
\newblock {\em On a type of matrix semigroups},
\newblock Semigroup Forum
\newblock {\bf  44} (1992), 125--128.

\bibitem{pon1}
---------,
\newblock {\em On matrix semigroups over a field ${K}$ conjugate to matrix
  semigroups over a proper subfield of ${K}$},
\newblock in Semigroups with Applications,  (J.~M. Howie, W.~D. Munn, and H.~J.
  Weinert, eds.),
\newblock World Sci., Singapore, 1992, pp.  1--5.

\bibitem{pra/sin}
D.~Prasad and K.~D. Singh,
\newblock {\em Matrix representation of semigroups},
\newblock J. Bihar Math. Soc.
\newblock {\bf  13} (1990), 45--48.

\bibitem{put1}
M.~S. Putcha,
\newblock {\em Matrix semigroups},
\newblock Proc. Amer. Math. Soc.
\newblock {\bf  88} (1983), 386--390.

\bibitem{putcha}
---------,
\newblock {\em Linear Algebraic Monoids},
\newblock Cambridge Univ. Press,
\newblock Cambridge, 1988.

\bibitem{rho/zal}
J.~Rhodes and Y.~Zalstein,
\newblock {\em Elementary representation and character theory of finite
  semigroups and its application},
\newblock in Monoids and Semigroups With Applications,  (J.~Rhodes, ed.),
\newblock World Sci., River Edge, 1991, pp.  334--367.

\bibitem{rog1}
A.~Rogers,
\newblock {\em A global theory of supermanifolds},
\newblock J.~Math. Phys.
\newblock {\bf  21} (1980), 1352--1365.

\bibitem{rog4}
---------,
\newblock {\em Graded manifolds, super\-manifolds and infinite-\-dimensional
  {G}rassmann algebras},
\newblock Comm. Math. Phys.
\newblock {\bf  105} (1986), 374--384.

\bibitem{rze1}
J.~Rzewuski,
\newblock {\em Supermatrix manifolds},
\newblock Rev. Math. Phys.
\newblock {\bf  29} (1991), 321--336.

\bibitem{schw1}
A.~S. Schwarz,
\newblock {\em To the definition of superspace},
\newblock Theor. Math. Phys.
\newblock {\bf  60} (1984), 37--42.

\bibitem{she}
I.~P Shestakov,
\newblock {\em Superalgebras as a building material for constructing
  counterexamples},
\newblock in Hadronic Mechanics and Nonpotential Interaction,
\newblock Nova Sci. Publ., Commack, NY, 1992, pp.  53--64.

\bibitem{siz}
W.~S. Sizer,
\newblock {\em Representations of semigroups of idempotents},
\newblock Czech. Math. J.
\newblock {\bf  30} (1980), 369--375.

\bibitem{urr/mor1}
L.~F. Urrutia and N.~Morales,
\newblock {\em The {C}ayley-{H}amilton theorem for supermatrices},
\newblock J.~Phys.
\newblock {\bf  {A27}} (1994), 1981--1997.

\bibitem{vla/vol}
V.~S. Vladimirov and I.~V. Volovich,
\newblock {\em Superanalysis. 1. {D}ifferential calculus},
\newblock Theor. Math. Phys.
\newblock {\bf  59} (1984), 3--27.

\bibitem{vol/aku}
D.~V. Volkov and V.~P. Akulov,
\newblock {\em On the possible universal neutrino interaction},
\newblock JETP Lett
\newblock {\bf  16} (1972), 621--624.

\bibitem{wes/zum1}
J.~Wess and B.~Zumino,
\newblock {\em Superspace formulation of supergravity},
\newblock Phys. Lett.
\newblock {\bf  B66} (1977), 361--364.

\bibitem{dewitt}
B.~S.~De Witt,
\newblock {\em Supermanifolds},
\newblock Cambridge Univ. Press,
\newblock Cambridge, 1984.

\bibitem{yan/bar}
S.~J. Yang and G.~P. Barker,
\newblock {\em Generalized {G}reen's relations},
\newblock Czech. Math. J.
\newblock {\bf  42} (1992), 211--224.

\bibitem{zal}
Y.~Zalstein,
\newblock {\em Studies in the representation theory of finite semigroups},
\newblock Trans. Amer. Math. Soc.
\newblock {\bf  161} (1971), 71--87.

\end{thebibliography}
\end{document}